\documentclass[12pt]{amsart}
\usepackage{amsthm, amssymb, amscd, mathptmx}

\DeclareMathAlphabet\mathscr{U}{rsfs}{m}{n}

\renewcommand{\Bbb}{\mathbb}

\newcommand{\cal}{\mathscr}

\newcommand{\C}{{\Bbb C}}

\newcommand{\PGL}{\mathop{\rm PGL}}
\newcommand{\GL}{\mathop{\rm GL}}

\renewcommand{\mod}{\mathop{\rm mod}}

\renewcommand{\O}{{\cal O}}

\renewcommand{\P}{{\Bbb P}}
\newcommand{\Pic}{\mathop{\rm Pic}}

\newcommand{\pt}{\mathop{\rm pt}}

\newcommand{\R}{{\Bbb R}}

\newcommand{\Z}{{\Bbb Z}}

\newcommand{\la}{\lambda}
\newcommand{\lra}{\longrightarrow}

\newcommand{\lma}{\longmapsto}
\newcommand{\p}{\prime}
\newcommand{\q}{\quad}

\renewcommand{\phi}{\varphi}

\newcommand{\eps}{\varepsilon}

\renewcommand{\theta}{\vartheta}

\newcommand{\ol}{\overline}

\theoremstyle{plain}
\newtheorem{theorem}{Theorem}[section]
\newtheorem{corollary}[theorem]{Corollary}
\newtheorem{proposition}[theorem]{Proposition}
\newtheorem{lemma}[theorem]{Lemma}

\theoremstyle{remark}
\newtheorem{remark}[theorem]{Remark}

\begin{document}

\title{Topological rigidity of algebraic $\P_3$-bundles over curves$^\dagger$}
\author{Alexander Schmitt}
\thanks{$^\dagger$ The paper is to appear in the Annali di Matematica Pura ed Applicata}
\address{Universit\"at GH Essen, FB6 Mathematik \& Informatik, D-45141 Essen, Germany.}
\email{\tt alexander.schmitt@uni-essen.de}

\maketitle

\begin{abstract}
A projective algebraic surface which is homeomorphic to a ruled
surface over a curve of genus $g\ge 1$ is itself a ruled surface
over a curve of genus $g$. In this note, we prove the analogous
result for projective algebraic manifolds of dimension $4$ in case
$g\ge 2$. 
\end{abstract}

\section{Introduction}
In the classification theory of projective algebraic manifolds, a very attractive
problem is the determination of all possible projective algebraic structures on a
given closed topological manifold. Among the classical results are the 
unicity of such structures on $\C\P_n$,
by Hirzebruch and Kodaira \cite{HK} for odd $n$ and Yau \cite{Y} in general,
and odd dimensional quadrics by Brieskorn \cite{Br1}.
Now, a class where one can hope to describe all projective algebraic structures
are projective bundles over curves. E.g., one knows that a surface which is 
homeomorphic to a ruled surface over a curve of genus $g\ge 1$ is itself
a ruled surface over a curve of the same genus $g$ (\cite{BPV}, (7.2), p.~202).
As far as the topological $4$-manifolds underlying rational surfaces are concerned, 
it is to the author's knowledge still an open problem
whether they might support surfaces of general type.
\par
In dimension three, Campana and Peternell \cite{CP} 
have shown that a projective algebraic
manifold $X$ which is homeomorphic to $\P_1\times\P_2$ is isomorphic either to
a manifold of the form $\P(\O_{\P_1}\oplus \O_{\P_1}(a)\oplus \O_{\P_1}(b))$ where 
$\pm(a+ b)\equiv 0\mod\ 3$ or to a manifold of the form $\P(E)$
where $E$ is a rank two vector bundle over $\P_2$ the Chern classes of which
satisfy $c_1(E)^2-4c_2(E)=0$. Freitag proves in \cite{Fr} that a manifold
which is homeomorphic to $\P(\O_{\P_1}^{\oplus 2}\oplus\O_{\P_1}(1))$, the blow up of
$\P_3$ in a line, is actually of the form
$\P(\O_{\P_1}\oplus \O_{\P_1}(a)\oplus \O_{\P_1}(b))$ with 
$\pm(a+ b)\equiv 1\mod\ 3$.
\par
In this note, we will derive the following result.
\begin{theorem}
\label{main}
Let $C$ be a smooth projective curve of genus $g\ge 2$ and $E$ a vector bundle
of rank $4$ and degree $k$ over $C$.
Suppose $X$ is a projective algebraic manifold which is homeomorphic to $\P(E)$.
Then, $X$ is of the form $\P(E^\p)$ where $E^\p$ is vector bundle of rank $4$
over a smooth projective curve $C^\p$ with $g(C^\p)=g$ and
$\deg E^\p\equiv \pm k \mod\ 4$.
\end{theorem}
Our strategy is quite similar to the one used in \cite{CP}.
First, we observe $\chi(\O_X)=\chi(\O_{\P(E)})=1-g$. The Miyaoka-Yau inequality thus
implies that $K_X$ can't be nef. Therefore, an elementary Mori contraction takes place.
The possible contractions are either those occuring in Ando's theorem (see Theorem~\ref{An} below)
or are contractions
onto a curve. The former contractions can be ruled out by topological arguments
while in the latter case a basic result from adjunction theory shows that the contraction
exhibits $X$ as a projective bundle over the given curve.
\par
A similar reasoning can be applied to describe those projective algebraic structures
on $\P_3$-bundles over $\P_1$ which are not of general type.
\begin{theorem}
\label{notsomain}
Suppose $X$ is a projective algebraic manifold which is homeomorphic to
$\P(\O_{\P_1}^{\oplus 3}\oplus \O_{\P_1}(k))$. Then, there are the following
possibilities:
\par
\begin{itemize}
\item[\rm i)] $X$ is of general type.
\par
\item[\rm ii)] $X$ is isomorphic to a manifold of the form
$\P(\O_{\P_1}\oplus \O_{\P_1}(a)\oplus \O_{\P_1}(b)
\oplus\O_{\P_1}(c))$ with $a+b+c\equiv \pm k\mod\ 4$.
\par
\item[\rm iii)] $k=0$, and 
$X$ is isomorphic to a manifold of the form $\P(E)$ where $E$ is a vector
bundle of rank two over $\P_3$ the Chern classes of which satisfy
$4c_2(E)-c_1(E)^2=0$.
\end{itemize}
\end{theorem}
\begin{remark}
In iii) we can normalize $E$, so that $c_1(E)=c_2(E)=0$.
Note that there are exactly two different topological $\C$ vector bundles of rank $2$
on $\P_3$ with trivial Chern classes, both of which carry holomorphic structures
\cite{OSS}. Call those two bundles $E$ and $E^\p$. Then,
$\P(E)$ and $\P(E^\p)$ are different as topological $\P_1$-bundles.
It would be interesting to see that they can't be abstractly homeomorphic.
This is not so easy to decide, because their standard topological invariants
(cohomology ring and characteristic classes) are equal. 
\end{remark}
One would expect the analog of Theorem~\ref{main} to hold also in the case $g=1$.
Unfortunately, I have not been able to establish it in that case.
Let me point out that, already for surfaces, the genus one case is harder 
(see \cite{BPV}).
In that case, an analysis of the Albanese map together with some specific
knowledge of fibrations of surfaces over curves is used. One might try the same
approach also in our situation, but it seems to me that the available information
might not be sufficient to arrive at the conclusion. Only if we assume that
$K_X$ be not nef, the result becomes easy.
\begin{proposition}
\label{notatallmain}
Let $C$ be a smooth projective curve of genus $g=1$ and $E$ a vector bundle
of rank $4$ and degree $k$ over $C$.
Suppose $X$ is a projective algebraic manifold which is homeomorphic to $\P(E)$
and for which $K_X$ is not nef.
Then, $X$ is of the form $\P(E^\p)$ where $E^\p$ is vector bundle of rank $4$
over a smooth projective curve $C^\p$ of genus one and
$\deg E^\p\equiv \pm k \mod\ 4$.
\end{proposition}
\section{Preliminaries}

\subsection*{$\P_3$-bundles over curves}
We recall some basic facts about $\P_3$-bundles over curves and adapt some observations
from \cite{Br} to our setting.
\par
Let $C$ be a smooth projective curve over $\C$ and $E$ a vector bundle of
rank $4$ and degree $k$. 
Let $\pi\colon\P(E)\lra C$ be the projective bundle of lines in $E$.
The Chern classes of $E$ are computed with the help of the exact sequence
$$
\begin{CD}
0 @>>> \O_{\P(E)} @>>> \pi^* E(1) @>>> T_{\P(E)} @>>> \pi^* T_C @>>> 0.
\end{CD}
$$
Set $s_E:= c_1(\O_{\P(E)}(1))$, and let $f_E$ be the cohomology class of a fibre of $\pi$.
Then,
\begin{eqnarray*}
c_1(E) &=& 4\cdot s_E+\bigl(k+2-2g(C)\bigr)\cdot f_E\\
c_2(E) &=& 6\cdot s_E^2 +\bigl(3k+8-8g(C)\bigr)\cdot s_Ef_E.
\end{eqnarray*}
\begin{lemma}
\label{coh}
Let $C_1$ and $C_2$ be two smooth projective curves and $E_1$ and $E_2$
two vector bundles of rank $4$ and degree $k_1$ and $k_2$ over $C_1$ and
$C_2$, respectively.
Then, the cohomology rings of $\P(E_1)$ and $\P(E_2)$ are isomorphic if and only
if $g(C_1)=g(C_2)$ and $k_1\equiv \pm k_2 \mod\ 4$.
\end{lemma}
\noindent
\begin{proof}
First of all, $b_1(\P(E_i))=b_1(C_i)=2g(C_i)$, $i=1,2$. Second, the classes
$a\cdot f_{E_i}$, $a\in\Z\setminus\{0\}$, are the only non trivial classes in $H^2(\P(E_i),\Z)$
with square zero, $i=1,2$. Thus, an isomorphism between $H^*(\P(E_1),\Z)$
and $H^*(\P(E_2),\Z)$ comes from an assignment
\begin{eqnarray*}
s_{E_1} &\lma & \pm s_{E_2} + l\cdot f_{E_2}
\\
f_{E_1} &\lma & \pm f_{E_2}
\end{eqnarray*}
for some $l\in\Z$. From this, one easily infers the claim.\qed
\end{proof}
\begin{proposition}
\label{reversy}
In the situation of Lemma~{\rm\ref{coh}}, the following holds true:
\begin{itemize}
\item[\rm i)]
$\P(E_1)$ and $\P(E_2)$ are orientation preservingly diffeomorphic,
if and only if $g(C_1)=g(C_2)$ and $k_1\equiv k_2 \mod\ 4$.
\item[\rm ii)]
$\P(E_1)$ and $\P(E_2)$ are diffeomorphic,
if and only if $g(C_1)=g(C_2)$ and $k_1\equiv \pm k_2 \mod\ 4$.
\end{itemize}
\end{proposition}
\noindent
\begin{proof}
i) follows immediately from the classification of topological
$\PGL_3(\C)$-bundles and the computations in the proof of Lemma~\ref{coh}.
\par
ii)
The necessity of the condition follows from Lemma~\ref{coh}. 
For its sufficiency, by i), it is enough to exhibit for every vector bundle
$E$ on a curve $C$ a projective bundle $\ol{P}$ which is orientation reversingly
diffeomorphic to $\P(E)$. To do so, let
$E$ be given w.r.t.\ to some covering $\{U_i\}_{i=1,...,t}$ of $C$
by the transition functions $\tau_{ij}\colon U_i\cap U_j\lra \GL_4(\C)$, $1\le i<j\le t$.
Let $\ol{E}$ be the differentiable $\C$ vector bundle defined by the complex conjugate
transition
functions $\ol{\tau}_{ij}$, $1\le i<j\le t$.
Then, the maps $f_i\colon U_i\times \P_3\lra U_i\times \P_3$, $\bigl(x,[z_0:\cdots:z_4]\bigr)
\lma \bigl(x, [\ol{z}_0:\cdots:\ol{z}_4]\bigr)$, $i=1,...,t$, glue to an orientation reversing
diffeomorphism $f\colon \P(E)\lra \P(\ol{E})$. 
\end{proof}
Finally, we note
\begin{lemma}
\label{auto}
Let $E$ be a rank $4$ vector bundle over the smooth projective curve $C$.
Then, $\P(E)$ possesses an orientation preserving self-diffeomorphism
which maps $s_E$ to $-s_E$ and $f_E$ to $-f_E$.
\end{lemma}
\noindent
\begin{proof}
We may assume that $C$ has a real structure with an $\R$-rational point $c_0\in C$
and that $E=\O_C^{\oplus 3}\oplus \O_C(k\cdot c_0)$.
Therefore, $\P(E)$ inherits a real structure, and complex conjugation on $\P(E)$
has the desired property. 
\end{proof}
\subsection*{The theorem of Ando}
This theorem from \cite{An} characterizes those elementary Mori contractions 
possessing
one dimensional fibres.
\begin{theorem}
\label{An}
Let $\rho\colon X \lra X^\p$ be an elementary Mori contraction from
a smooth variety $X$ of dimension $n$. Let $F=\rho^{-1}(z)$ be
a one dimensional fibre.
\begin{enumerate}
\item[\rm i)] If $\rho$ is birational, then $F$ is isomorphic to $\P_1$ and
$-K_X.[F]=1$. In this case,
$X^\p$ is smooth, and $\rho$ is the blow up along a smooth 
subvariety of codimension 2.
\item[\rm ii)] If $\rho$ is not birational, then
$X^\p$ is smooth of dimension $n-1$,
and $\rho$ is a flat conic bundle. In particular, any smooth curve $F$
which is contracted by $\rho$ is isomorphic to $\P_1$, and 
$-K_X.[F]\in \{\, 1,2\,\}$. Here, $-K_X.[F]=1$ happens if and only if $F$ is contained
in a singular fibre, i.e., either a reducible fibre or a fibre of
multiplicity two.
\end{enumerate}
\end{theorem}
\subsection*{A result from adjunction theory}
The following theorem gives a handy criterion for a manifold to be a projective
bundle over a curve.
\begin{theorem}
\label{Adj}
Let $X$ be a projective algebraic manifold of dimension $n\ge 3$ and $L$ an ample
line bundle on $X$. Suppose $b_2(X)\ge 2$. If $K_X+(n-1)\cdot L$ is not nef, then
$X$ is a projective bundle over a smooth curve $C$ and $L=\O_X(1)$.
\end{theorem}
\noindent
\begin{proof} See \cite{Fu}, Theorem~(11.7), p.\ 96. \end{proof}
\subsection*{A characterization of $\P_3$}
The result below is taken from the paper \cite{LS}.
\begin{theorem}
\label{LS}
Let $X$ be a projective algebraic manifold of dimension $3$.
Assume that the cohomology ring of $X$ is isomorphic to the one
of $\P_3$. Then, $X$ is isomorphic to $\P_3$ as algebraic manifold.
\end{theorem}
\section{Proof of Theorem~\ref{main}}
By Proposition~\ref{reversy}~ii), we may assume that we have an orientation preserving
homeomorphism $h\colon X\lra \P(E)$. Let $s$ and $f$ be the images of $s_E$ and $f_E$
under the isomorphism $h^*\colon H^2(\P(E),\Z)\lra H^2(X,\Z)$.
The Betti numbers of $X$ are
$$
b_0(X)=b_8(X)=1,\q b_{2i}(X)=2,\ i=1,2,3,\q b_1(X)=b_7(X)=2g,
$$
$$
b_3(X)=b_5(X)=4g.
$$
Thus, $h^1(\O_X)=g$, $h^2(\O_X)=h^4(\O_X)=0$, in particular, $s$ and $f$ belong
to $H^{1,1}(X)$. Since $H^3(X,\Z)= H^1(X,\Z)\cdot s\oplus H^1(X,\Z)\cdot f$, we find
$h^{2,1}(\O_X)=h^{1,2}(\O_X)=2g$, i.e., $h^3(\O_X)=0$.
We conclude
$$
\chi(\O_X)\q=\q\chi(\O_{\P(E)})\q=\q 1-g.
$$
Let $c_i$, $i=1,2,3,4$, be the Chern classes of $X$ and
$p_1=c_1^2-2c_2$ and $p_2=c_2^2-2c_1c_3+2c_4$ its Pontrjagin classes.
By Riemann-Roch, we have
\begin{eqnarray*}
\chi(\O_X) &=& -\frac{1}{720}\bigl(c_1^4-4c_1^2c_2-3c_2^2-c_1c_3+c_4\bigr)
\\
&=& -\frac{1}{720}\left(p_1^2+\frac{1}{2}p_2-\frac{15}{2}c_2^2\right).
\end{eqnarray*}
Since the Pontrjagin numbers
are invariants of the oriented homeomorphy type~\cite{Nov},
we infer
\begin{equation}
\label{inv1}
c_2^2\q=\q c_2(\P(E))^2\q=\q 96(1-g)
\end{equation}
and
\begin{equation}
\label{inv2}
4c_1^2c_2-c_1^4\q=\q 4c_1(\P(E))^2c_2(\P(E))-c_1(\P(E))^4\q=\q 384(1-g).
\end{equation}
This has two important consequences:
\begin{corollary}
\label{MY}
The canonical bundle $K_X$ of $X$ cannot be nef.
\end{corollary}
\noindent
\begin{proof}
If $K_X$ were nef, then the Miyaoka-Yau inequality \cite{MY} would yield
$3c_1^2c_2-c_1^4\ge 0$. Thus, also $4c_1^2c_2-c_1^4\ge 0$.
But we have seen in (\ref{inv2}) that this quantity is negative. 
\end{proof}
Now, let $c_1=a\cdot s+b\cdot f$. By Lemma~\ref{auto}, we may assume $a\ge 0$.
\begin{corollary}
\label{Est}
We have\q
$
a\ge 4\q\hbox{and}\q 4b-ak\q <\q 0.
$
\end{corollary}
\noindent
\begin{proof}
By the topological invariance of the Pontrjagin classes~\cite{Nov}, 
$p_1=4\cdot s+ 2k\cdot f$. Thus,
$$
c_2\q=\q\frac{1}{2}\bigl(a^2-4\bigr)\cdot s^2+\bigl(ab-k\bigr)\cdot sf
$$
and
$$
c_2^2\q=\q a(a^2-4)\left(b-\frac{ak}{4}\right)\q\stackrel{(\ref{inv1})}{=}\q 96(1-g).
$$
Observe that $a$ must be an even number, so the claim follows. 
\end{proof}
Since $K_X$ is not nef, there exists an elementary Mori contraction
$\rho\colon X\lra X^\p$.
\begin{lemma}
\label{cont}
If $\dim X^\p\ge 2$, then $\rho$ doesn't contract any surface.
\end{lemma}
\noindent
\begin{proof}
Let $L$ be an ample line bundle on $X^\p$ and $c_1(\rho^* L)=\alpha\cdot s+\beta\cdot f$.
Assume $F$ is an irreducible surface with $\rho(F)=\{\pt\}$. Write
$[F]=\gamma\cdot s^2+\delta\cdot sf$. We find
$$
0\q=\q (\alpha\cdot s+\beta\cdot f).(\gamma\cdot s^2+\delta\cdot sf)
\q=\q \alpha\gamma\cdot s^3+(\alpha\delta+\beta\gamma)\cdot s^2f.
$$ 
Thus, $\alpha\gamma=0$. Since we assume $\dim X^\p\ge 2$, $\alpha\neq 0$, whence
$\gamma=0$ and $\delta=0$, a contradiction. 
\end{proof}
We show first that $\rho$ cannot be birational. If $\rho$ were birational,
then, by Lemma~\ref{cont} and Ando's theorem, $\rho$ would be the blow
up of a smooth $4$-manifold $X^\p$ in a smooth surface.
Let $L$ be an ample line bundle on $X^\p$.
Denote by $E$ the exceptional divisor of $\rho$. With $\la:=
\alpha\cdot s+\beta\cdot f:=\pm c_1(\rho^* L)$ and
$\eps:=\gamma\cdot s+\delta\cdot f:=\pm [E]$, we find a basis for $H^2(X,\Z)$,
such that 
\begin{equation}
\label{exz1}
\alpha\delta-\beta\gamma=1
\end{equation}
and
$$
0\q=\q\la^3.\eps\q=\q(\alpha\cdot s+\beta\cdot f)^3.(\gamma\cdot s+\delta\cdot f)
\q=\q \alpha^2(-\alpha\gamma k+\alpha\delta+3\beta\gamma).
$$
Again $\alpha\neq 0$, so that
\begin{equation}
\label{exz2}
-\alpha\gamma k+\alpha\delta+3\beta\gamma\q=\q 0.
\end{equation}
From (\ref{exz1}) and (\ref{exz2}), one infers
$$
\gamma(-\alpha k+4\beta)\q=\q-1\qquad {\rm and}\qquad -\alpha(\gamma k-4\delta)\q=\q 3.
$$
Hence, $\gamma=\pm 1$ and $\alpha\in\{\pm 1,\pm 3\}$.
Without loss of generality, we may assume $\gamma=1$.
Then, $\beta=(\alpha k-1)/4$.
Observe that for even $k$, we have $w_2(X)=0$, i.e., $c_1(X)$ is divisible
by two. Therefore, no curve $F$ with $c_1(X).F=1$ can exist. Since we can clearly assume
$k\in \{\, 0,1,2,3\,\}$, we only have to deal with the cases $k=1$ and $k=3$.
For $k=1$, the possible solutions are $(\alpha,\beta)=(-3,-1)$ and 
$(\alpha,\beta)=(1,0)$. In the latter case, $L^4=\la^4=s^4=-1$ which is impossible.
In the former case, we find $\eps=s$.
Therefore, for some $l>0$,
$$
c_1(X)\q=\q (3l\pm 1)\cdot s + l\cdot f.
$$
Since $(4l-3l\mp 1)\ge 0$, this contradicts Corollary~\ref{Est}.
\par
For $k=3$, $(\alpha,\beta)=(-1,-1)$ and $(\alpha,\beta)=(3,2)$ are the possible
solutions. The second solution is excluded as 
$L^4=\la^4=-27$. In the first case, $\eps=s$, so that, 
for some $l\ge 3$,
$$
c_1(X)\q=\q (l\pm 1)\cdot s + l\cdot f.
$$
This is again not an option.
\par
Assume now $\dim X^\p=3$. Then, $X^\p$ is smooth, and $\rho$
is a conic bundle. 
Let $L$ be an ample line bundle on $X^\p$ and
$\alpha\cdot s+\beta\cdot f:=c_1(\rho^* L)$.
Then,
$$
0\q=\q (\alpha\cdot s+\beta\cdot f)^4\q=\q \alpha^3(-\alpha k+ 4\beta)
$$
and $\alpha\neq 0$.
Let $F\cong \P_1$ be a smooth curve which is contracted by $\rho$.
One easily computes
that the class $[F]$ is a multiple of $s^3+(3/4)k\cdot s^2f$.
Therefore, the composite map
$$
\P_1\cong F\hookrightarrow X \stackrel{h}{\lra} \P(E) \lra C
$$
has non zero degree, an absurdity (cf.~\cite{BPV}, Cor.~(1.2), p.~11).
\par
We finally arrive at the case where $X^\p=C^\p$ is a smooth projective curve.
Then, for any point $c\in C^\p$, $\rho^*(\O_C(c))=\alpha\cdot f$, for some $\alpha\in\Z
\setminus\{0\}$.
Since $-K_X$ is $\rho$-ample, we see
$$
0\q<\q (a\cdot s+b\cdot f)^3.(\alpha\cdot f)\q=\q a^3\alpha, 
$$
i.e., $\alpha>0$. For $l\gg 0$, the line bundle
$$
L\q:=\q s+l\cdot f\q=\q
\frac{1}{a} c_1(X)+ \left(l-\frac{b}{a}\right)\cdot f
$$
is ample. The bundle $K_X+aL=(al-b)\cdot f$ is nef,
but $H:=K_X+(a-1)L=-s+((a-1)l-b)\cdot f$ isn't, because
$$
H^3.f\q=\q -s^3f\q=\q -1 \q <\q 0.
$$
Thus, since $a\ge 4$, by Corollary~\ref{Est}, we must actually have $a=4$, and we can conclude
by Theorem~\ref{Adj}.\qed
\section{Proof of Theorem~\ref{notsomain}}
We may assume $k\in\{\, 0,1,2,3\,\}$.
Let $c_1(X)=a\cdot s+b\cdot f$ with $a\ge 0$.
\begin{lemma}
\label{Chern}
Either $a=4$ and $b=2+k$ or $k=3$ and $a=6$ and $b=5$.
\end{lemma}
\noindent
\begin{proof}
As before, we find
$$
a(a^2-4)\left(b-\frac{ak}{4}\right)\q=\q 96.
$$
Now, $a=2\alpha$ is an even number. If $k$ is even, $b^\p:= b-(ak/4)$ is an integer, and
we have
$\alpha(\alpha^2-1)b^\p=12$ which forces $\alpha=2$, whence $a=4$ and $b=2+k$ as desired.
If $k$ is odd, then $b^\p:=2b-(ak/2)$ is an integer, and
$\alpha(\alpha^2-1)b^\p=24$ can be solved by $\alpha=2$, leading to $a=4$ and $b=2+k$
as before, or by $\alpha=3$, leading to $a=6$ and $b=(1/2)+(3k/2)$.
For $k=1$, this gives $c_1(X)=6\cdot s+ 2\cdot f$, i.e., $w_2(X)=0$, a contradiction.
 
\end{proof}

Since 
\begin{eqnarray*}
(4\cdot s+ (2+k)\cdot f)^4 & =& 512\q>\q 0\q \hbox{and}
\\
(6\cdot s+  5\cdot f)^4 &=& 432 \q>\q 0,
\end{eqnarray*}
$X$ is either of general type, or $K_X$ cannot be nef.
We assume the latter case.
\par
Thus, we have again an elementary Mori contraction $\rho\colon X\lra X^\p$, and
we use the same notation as before.
Note that Lemma~\ref{cont} and most of the computations afterwards still apply.
Assume first that $\rho$ is birational. Then, $k=1$ or $k=3$, because otherwise
$w_2(X)=0$.
For $k=1$, we have, by Lemma~\ref{Chern},
$c_1(X)=4\cdot s+ 3\cdot f$.
Let $L$ be an ample line bundle on $X^\p$. We have already seen that $\rho^*(c_1(L))$
is a multiple of $3\cdot s+f$. Let $F$ be a fibre
by $\rho$. For this curve, $(3\cdot s+f).[F]=0$ and $(4\cdot s+3\cdot f).[F]=1$, by 
Theorem~\ref{An}.
But one verifies that this is impossible.
\par
Now, suppose $k=3$ and
let $\sigma$ be the class with $\rho^*\sigma=s+f$. Since
$K_X=\rho^*(K_{X^\p})+E$ (\cite{Ful}, Ex.~15.4.3), we find $[K_{X^\p}]=-5\sigma$. If $X$ is not of general
type, then the same holds for $X^\p$. It follows that $X^\p\cong\P_4$ and
$\sigma=c_1(\O_{\P_4}(1))$ (\cite{Fu}, Theorem~(11.2), p.~93).
Next, with $[S]$ the class of the surface in which $X$ is blown up, 
\begin{eqnarray*}
c_2(X) &\stackrel{\hbox{\cite{Ful}, (15.4.3)}}{=} & \rho^*(c_2(X^\p)+[S])-\rho^*c_1(X^\p).E
\\
       &= & 10(s+f)^2\pm(5s+5f).s + \rho^*[S].
\end{eqnarray*}
Here, we have the "$+$"-sign for $c_1(X)=4\cdot s+5\cdot f$ and the
"$-$"-sign for $c_1(X)=6\cdot s+5\cdot f$. It follows in both cases that
$[S]=\sigma^2$, i.e., $X$ is the blow-up of $\P_4$ in a linearly embedded
$\P_2$. The pencil of hyperplanes passing through that $\P_2$ induces on $X$
the structure of a projective bundle over $\P_1$, and we are done in this case.
\par
Now, we turn to the case where $X$ is a conic bundle over a smooth $3$-manifold.
We first look at the case $c_1(X)=4\cdot s+(2+k)\cdot f$ and 
claim that $k=0$, i.e, $X$ is homeomorphic to $\P_1\times\P_3$.
As usual, $L$ is an ample generator for $H^2(X^\p,\Z)$
and $\alpha\cdot s+\beta\cdot f:=c_1(\rho^* L)$.
As before, we compute that, for a smooth curve $F$
which is contracted by $\rho$, $[F]=l(s^3+(3/4)k\cdot s^2f)$ for some $l\in\Z$.
We choose $F$ with $c_1(X).[F]=2$. As
$$
l(4\cdot s + (2+k)\cdot f).\left(s^3+\frac{3}{4}k\cdot s^2f\right)\q=\q 2l,
$$ 
$l$ must be equal to $1$. 
But $s^3+(3/4)k\cdot s^2f$ is an integral class if and only if $k=0$.
In this case, all fibres must be smooth, i.e., $X$ is a projective bundle over $X^\p$,
because $H^2(X^\p,\O_{X^\p}^*)=H^3(X^\p,\Z)=0$.
Since $\rho^*\colon H^*(X^\p,\Z)\lra H^*(X,\Z)$ is injective, it follows readily that
$$
H^*(X^\p,\Z)\q\cong\q \Z[s]/\langle s^4\rangle.
$$
By Theorem~\ref{LS}, $X^\p\cong \P_3$. We see that $X\cong \P(E)$ where $E$ is
a vector bundle of rank $2$ over $\P_3$. It is an easy exercise to check
that $4c_2(E)-c_1(E)^2=0$.
\par
Now, we assume $c_1(X)=6\cdot s+ 5\cdot f$ and $k=3$. The same argument as
before shows that the class $[F]$ of a smooth curve contracted by $\rho$
is a multiple of $s^3+(9/4)\cdot s^2f$. As
$$
(6\cdot s+ 5\cdot f).\left(s^3+\frac{9}{4} \cdot s^2f\right)\q=\q \frac{1}{2},
$$
we conclude again that $\rho\colon X\lra X^\p$ is a projective bundle,
and that the cohomology class of a fibre is given by $4\cdot s^2+9\cdot s^2f$.
Observe that $\chi(\O_{X^\p})=1$. This implies that $X^\p$ is a Fano manifold.
We obviously have $\Pic X^\p=H^2(X^\p,\Z)\cong \Z$ and $b_3(X^\p)=0$.
Let $L$ be an ample generator of $\Pic X^\p$. We know already that
$c_1(\rho^*L)$ is a multiple of $s+ (3/4)\cdot f$. Since $c_1(\rho^* L)$
forms part of a $\Z$-basis for $H^2(X,\Z)$, this class must equal to $4\cdot s+ 3\cdot f$.
In particular,
$$
c_1(\rho^* L)^3\q=\q 64\cdot s^3 + 144\cdot s^2f.
$$
Comparing this with the formula for the class $[F]$ of a fibre, we find
$L^3=16$. The classification of Fano $3$-manifolds \cite{IP}
shows that no Fano $3$-manifold
of degree $16$ with $b_3=0$ exists, so that this case cannot occur.
\par
The remaining case of a contraction onto a curve, this time $\P_1$, is handled
as before.\qed
\section{Proof of Proposition~\ref{notatallmain}}
With the same methods as before, one first establishes
\begin{lemma}
\label{Est2}
For $c_1(X)=a\cdot s+b\cdot f$, one has either $a=0$ or $a\ge 4$ and $b=ak/4$.
\end{lemma}
From the Betti numbers of $X$, we read off that $h^{1,0}(X)=1$, i.e.,
the Albanese torus of $X$ is an elliptic curve. Since the image of the
Albanese map generates the Albanese torus as a group, the Albanese 
map ${\rm Alb}\colon X\lra {\rm Alb}(X)$ is surjective. Let $h\colon X\lra A$
be the Stein factorization of ${\rm Alb}$ with $A$ a smooth connected curve.
Since $K_X$ is not nef, there are smooth rational curves $F$ with $K_X.[F]<0$. These must
be contracted by $h$, because they are contracted by ${\rm Alb}$.
As we have observed before, the integer multiples of $f$ are the only classes
with square zero. This implies that the classes of the fibres of $h$ are multiples
of $f$ and that the class of any curve which is contracted by $h$ is a multiple
of $s^2.f$. In particular, $[F]=c\cdot s^2.f$ for some $c\in\Z\setminus\{0\}$.
We conclude that, in Lemma~\ref{Est2}, we must have $a\ge 4$ and that $-K_X$
is $h$-ample. From here on, we may proceed as before. \qed

\end{document}